\documentclass[prl, showpacs, twocolumn, superscriptaddress, amsmath, amssymb, floatfix]{revtex4}

\usepackage{graphicx}
\usepackage{array}
\usepackage{amssymb}
\usepackage{amsfonts}
\usepackage{amsmath}
\usepackage{mathrsfs}
\usepackage{color}
\usepackage{booktabs}
\usepackage{threeparttable}
\usepackage{multirow}
\usepackage{subfigure}
\usepackage{times}
\usepackage{epsfig}
\usepackage{threeparttable}
\usepackage{chngpage}
\usepackage{float}
\usepackage{color}

\begin{document}

\title{Individual dynamics induces symmetry in network controllability}

\author{Chen Zhao}
\affiliation{School of Systems Science, Beijing Normal University,
Beijing, 10085, P. R. China}

\author{Wen-Xu Wang}
\affiliation{School of Systems Science, Beijing Normal University,
Beijing, 10085, P. R. China}

\author{Yang-Yu Liu}\email{yyl@channing.harvard.edu}
\affiliation{Channing Division of Network Medicine, Brigham and Women's
Hospital, Harvard Medical School, Boston, Massachusetts 02115, USA}
\affiliation{Center for Complex Network Research and Department of
Physics,Northeastern University, Boston, Massachusetts 02115, USA}

\author{Jean-Jacques Slotine}\email{jjs@mit.edu}
\affiliation{Nonlinear Systems Laboratory, Massachusetts Institute of
Technology, Cambridge, Massachusetts, 02139, USA}
\affiliation{Department of Mechanical Engineering and Department of Brain and
Cognitive Sciences, Massachusetts Institute of Technology, Cambridge, Massachusetts, 02139, USA}

\begin{abstract}
Controlling complex networked systems to a desired state is a key research goal in contemporary science.
Despite recent advances in studying the impact of network
topology on controllability, a comprehensive understanding of the synergistic effect of
network topology and individual dynamics on controllability
is still lacking.
Here we offer a theoretical study with particular interest in the
diversity of dynamic units characterized by different types of individual dynamics.
Interestingly, we find a global symmetry accounting for the invariance of controllability
with respect to exchanging the densities of any two different types of dynamic units, irrespective of
the network topology. The highest controllability arises at the global symmetry
point, at which different types of dynamic units are of the same
density.
The lowest controllability occurs when all self-loops are either
completely absent or present with identical weights.
These findings further improve our understanding of
network controllability and have implications for devising
the optimal control of complex networked systems in a wide range of fields.
\end{abstract}


\maketitle


As a key notion in control theory, controllability denotes our
ability to drive a dynamic system from any initial state to any
desired final state in finite
time~\cite{Kalman-JSIAM-63,Luenberger-Book-79}.
For the canonical linear time-invariant (LTI) system $\dot{\bf x}= A
{\bf x} + B {\bf u}$ with state vector ${\bf x} \in
\mathbb{R}^N$, state matrix $A \in \mathbb{R}^{N\times N}$ and control
matrix $B \in \mathbb{R}^{N\times M}$,  Kalman's rank
condition $\text{rank}[B,AB,\cdots, A^{N-1}B] = N$ is sufficient and
necessary to assure controllability.
Yet, in many cases system parameters are not exactly known, rendering
classical controllability tests impossible.
By assuming that system parameters are either fixed zeros or freely
independent,  structural control theory (SCT) helps us overcome this
difficulty for linear time-invariant systems~\cite{Lin:1974,Shields-IEEE-76,Hosoe-IEEE-80,
  SContrl1,Dion-Automatica-03}.
Quite recently, many research activities have been devoted to study the
structural controllability of systems with complex network structure,
where system parameters (e.g., the elements in $A$, representing link
weights or interaction strengths between nodes) are typically not
precisely known, only the zero-nonzero pattern of $A$ is known~\cite{LSB:2011,NV:2012,LSB:pOne,WNLG:2012,
  Kurth:2012,WGG:2012, PLSB:2013,bimodality:2013}.
Network controllability problem can be typically posed as a combinatorial
optimization problem, i.e., identify a minimum set of driver nodes,
with size denoted by $N_\mathrm{D}$, whose control is sufficient to
fully control the system’s dynamics~\cite{LSB:2011}.
Other controllability related issues, e.g., energy cost, have also
been extensively studied for complex networked
systems~\cite{Yan:2012,GSZPB:2012, LSB:observability,SM:PRL}.
While the intrinsic individual dynamics can be incorporated in the
network model, it would be more natural and fruitful to consider their
effect separately.
Hence, most of the previous studies focused on the
impact of network topology, rather than the individual dynamics of
nodes, on network controllability~\cite{LSB:2011,WNLG:2012}.

If one explores the impact of individual dynamics on network
controllability in the SCT framework, a specious result would be
obtained --- a single control input can make an arbitrarily large linear
system controllable. Although this result as a
special case of the minimum inputs theorem can be proved~\cite{LSB:2011} and its implication was further emphasized in~\cite{Cowan:2012},
this result is inconsistent with empirical situations, implying that the SCT is inapplicable
in studying network controllability, if
individual dynamics of nodes are imperative to be incorporated to capture the collective dynamic behavior of a networked system.
To overcome this difficulty, and more importantly, to understand the
impact of individual dynamics on network controllability, we revisit the key assumption of SCT, i.e.,
the independency of system parameters. We anticipate that major new
insights can be obtained by relaxing this assumption, e.g.,
considering the natural diversity and similarity of individual dynamics. This
also offers a more realistic characterization of many real-world
networked systems where not all the system parameters are completely
independent.

To solve the network controllability problem with dependent system
parameters, we rely on the recently developed exact
controllability theory (ECT)~\cite{exact_control}. ECT enables us to systematically explore the role of individual
dynamics in controlling complex systems with arbitrary network
topology.
In particular, we consider prototypical linear forms of individual
dynamics (from first-order
to high-orders) that can be incorporated within
the network representation of the whole system
in a unified matrix form.
%
%
This paradigm leads to the discovery of a striking symmetry in network controllability: if we exchange the
fractions of any two types of dynamic units, the system's
controllability (quantified by $N_\mathrm{D}$) remains the same.
This exchange-invariant property gives rise to a global symmetry
point, at which the highest controllability (i.e., lowest number of
driver nodes) emerges. This symmetry-induced optimal controllability
holds for any network topology and various categories of individual dynamics.
We substantiate these findings numerically in a variety of network models.\\

\begin{figure}[t!]
\begin{center}
\epsfig{figure=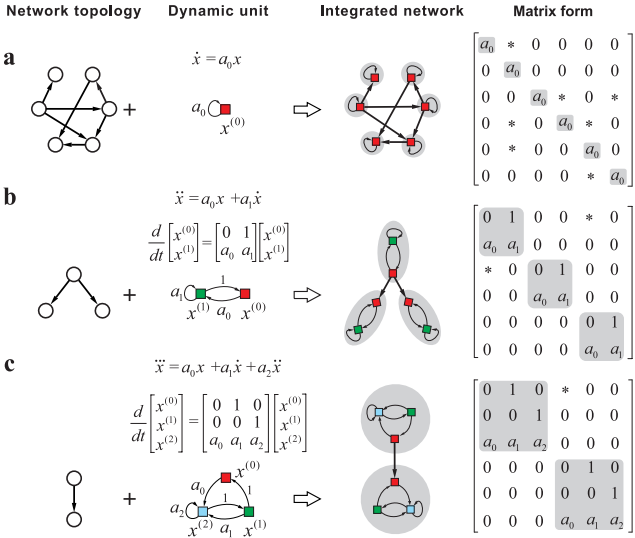,width=\linewidth}
\caption{Integration of network topology and (a)
1st-order, (b) 2nd-order and (c) 3rd-order intrinsic individual dynamics. For a $d$th-order individual dynamics
$x^{(d)}=a_0x^{(0)}+a_1x^{(1)}+\cdots + a_{d-1}x^{(d-1)}$,
we denote each order by a colored square and the couplings among orders
are characterized by links or self-loops. This graphical representation
allows individual dynamics to be integrated with their coupling network topology,
giving rise to a unified matrix that reflects the dynamics of
the whole system. In particular, each dynamic unit in the unified matrix
corresponds to a diagonal block and the nonzero elements (denoted by $*$)
apart from the blocks stand for the couplings among
different dynamic units. Therefore, the original network consisting of $N$ nodes with
order $d$ is represented in a $dN \times dN$ matrix.
}
\label{fig:illus_motif}
\end{center}
\end{figure}

Exact controllability theory (ECT)~\cite{exact_control} claims that for arbitrary network topology
and link weights characterized by the state matrix $A$
in the LTI system $\dot{\mathbf{x}}= A\mathbf{x}+B\mathbf{u} $, the minimum number of
driver nodes $N_\text{D}$ required to be controlled by imposing independent signals
to fully control the system is given by the maximum geometric multiplicity
$\max_i\{\mu(\lambda_i)\}$ of $A$'s eigenvalues $\{\lambda_i\}$~\cite{multi1,multi2,multi3,multi4,multi5}. Here $\mu(\lambda_i) \equiv N-
\text{rank}(\lambda_i I_N-A)$ is the geometric multiplicity of the eigenvalue $\lambda_i$ and $I_N$ is the identity matrix.
Calculating all the eigenvalues of $A$ and
subsequently counting their geometric multiplicities are generally
applicable but computationally prohibitive for large networks. If $A$ is symmetric, e.g., in undirected
networks, $N_\text{D}$ is simply given by the maximum algebraic multiplicity $\max_i\{\delta(\lambda_i)\}$,
where $\delta(\lambda_i)$ denotes the degeneracy of eigenvalue $\lambda_i$. Calculating $N_\text{D}$
in the case of symmetric $A$ is more computationally affordable than in
the asymmetric case. Note that for structured systems where the elements
in $A$ are either fixed zeros or free independent parameters, ECT offers the
same results as that of the SCT~\cite{exact_control}. \\

We first study the simplest case of first-order individual dynamics $\dot{x}_i=a_0x_i$.
The dynamical equations of a linear time-invariant control system associated with first-order individual
dynamics~\cite{Slotine:book} can be written as
\begin{equation} \label{eq:1st_order_node}
\dot{\mathbf{x}}=\Lambda\mathbf{x}+ A\mathbf{x}+B\mathbf{u} = \Phi\mathbf{x}+B\mathbf{u},
\end{equation}
where the vector $\mathbf{x} = (x_1,\cdots,x_N)^\text{T}$ captures the states of $N$ nodes,
$\Lambda\in \mathbb{R}^{N\times N}$ is a diagonal matrix representing intrinsic
individual dynamics of each node,
$A\in \mathbb{R}^{N\times N}$ denotes the coupling matrix or the weighted
wiring diagram of the networked system, in which $a_{ij}$ represents
the weight of a directed link from node $j$ to $i$ (for undirected networks, $a_{ij}=a_{ji}$).
$\mathbf{u}=(u_1,u_2,\cdots,u_M)^\text{T}$ is the input vector of $M$ independent signals,
$B \in \mathbb{R}^{N\times M}$ is the control matrix, and $\Phi \equiv \Lambda + A$
is the state matrix. Without loss of generality, we assume $\Lambda$ is
a ``constant" matrix over the field $\mathbb{Q}$ (rational numbers), and $A$ is a
structured matrix over the field $\mathbb{R}$ (real numbers). In other words, we
assume all the entries in $\Phi$ have been rescaled by the individual dynamics parameters.
The resulting state matrix $\Phi$ is usually called a \emph{mixed matrix} with respect to
$(\mathbb{Q}, \mathbb{R})$~\cite{mixed_matrix}. The first-order
individual dynamics in $\Phi$ is captured by self-loops in the network
representation of $\Phi$ (see Fig.~\ref{fig:illus_motif}a).
$N_\text{D}$ can then be determined by calculating
the maximum geometric multiplicity $\max_i\{ \mu(\lambda_i)\}$
of $\Phi$'s eigenvalues.

\begin{figure}[t!]
\begin{center}
\epsfig{figure=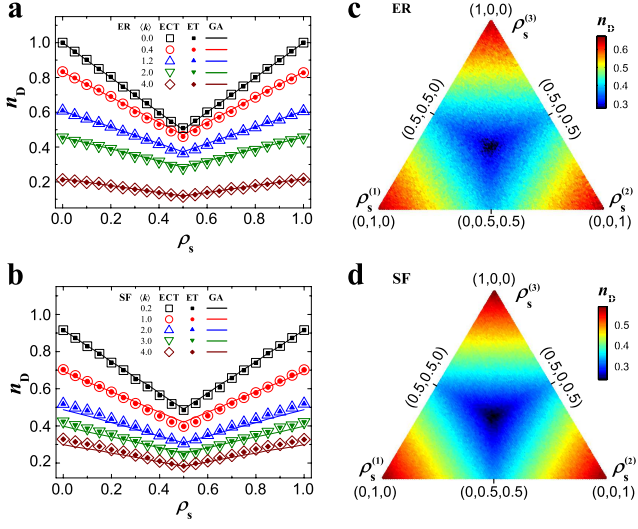,width=\linewidth}
\caption{(a)-(b) controllability measure $n_\text{D}$
in the presence of a single type of nonzero self-loops with fraction $\rho_\text{s}$ for random (ER)
networks (a) and scale-free (SF) networks (b) with different average degree $\langle k\rangle$.
(c)-(d) $n_\text{D}$ of ER (c) and SF networks ({\bf d}) with three types of self-loops
$s_1$, $s_2$ and $s_3$ with density $\rho_\text{s}^{(1)}$, $\rho_\text{s}^{(2)}$ and $\rho_\text{s}^{(3)}$,
respectively. ECT denotes the results obtained from the exact
controllability theory, ET denotes the results obtained from the efficient tool and GA denotes
the results obtained from the graphical approach (see SM Sec.3).
The color bar denotes the value of $n_\text{D}$ and the coordinates in the triangle stands
for $\rho_\text{s}^{(1)}$ $\rho_\text{s}^{(2)}$ and $\rho_\text{s}^{(3)}$. The networks are described by structured matrix $A$ and their
sizes in (a)-(d) are 2000. The results from ECT and ET are averaged over 30 different realizations, and those from GA are over 200 realizations.}
\label{fig:1st_order}
\end{center}
\end{figure}

We study two canonical network models (Erd\"os-R\'enyi and Scale-free) with
random edge weights and a $\rho_s$ fraction of nodes associated with identical
individual dynamics (i.e., self-loops of identical weights).
As shown in Fig.~\ref{fig:1st_order}a,b,
the fraction of driver nodes $n_\text{D}\equiv N_\text{D}/N$
is symmetric about $\rho_\text{s}=0.5$, regardless of the network topology.
Note that the symmetry cannot be predicted by SCT in the sense
that in case of completely independent
self-loop weights $n_\text{D}$ will monotonically decrease to $1/N$ as
$\rho_\text{s}$ increases to 1, implying that a single driver
node can fully control the whole network~\cite{Cowan:2012}.
The symmetry can be theoretically predicted (see SM Sec.2.2).
An immediate but counterintuitive result from the symmetry is
that $n_\text{D}$ in the absence of self-loops is exactly the same as the case that each node has a
self-loop with identical weight. This is a direct consequence of
Kalman's rank condition for controllability~\cite{Kalman-JSIAM-63}:
\begin{widetext}
\begin{equation}
\text{rank}[B,AB,\cdots, A^{N-1}B] =
\text{rank}[B,(A+w_\text{s}I_N)B,\cdots, (A+w_\text{s}I_N)^{N-1}B]
\label{eq:kalman_equal}
\end{equation}
\end{widetext}
where the left and the right hand sides are the rank of controllability matrix in the absence and full
of identical self-loops, respectively (see SM Sec.1 for proof).

The presence of two types of nonzero self-loops $s_2$ and $s_3$ leads to even richer behavior of controllability.
If the three types of self-loops (including self-loops of zero weights) are randomly
distributed at nodes, the impact of their fractions on $n_\text{D}$ can be
visualized by mapping the three fractions into a 2D triangle (or
2-simplex), as shown in Fig.~\ref{fig:1st_order}c,d. We see that $n_\text{D}$ exhibits symmetry in the
triangle and the minimum $n_\text{D}$ occurs at the center that represents
identical fractions of the three different self-loop types. The symmetry-induced
highest controllability can be generalized to arbitrary number of self-loops.
Assume there exist $n$ types of self-loops $s_1,\cdots ,s_n$ with weights
$w_\text{s}^{(1)},\cdots ,w_\text{s}^{(n)}$, respectively, we have
\begin{eqnarray}
N_\text{D}= N- \min_i \bigg\{ \text{rank}\big(\Phi-w_\text{s}^{(i)} I_N \big) \bigg\}
\label{eq:unify_1st}
\end{eqnarray}
for sparse networks with random weights
(see SM Sec. 2 for detailed derivation and the formula of dense networks).
An immediate prediction of Eq.~(\ref{eq:unify_1st}) is that $N_\text{D}$
is primarily determined by the self-loop with the highest density,
simplifying Eq.~(\ref{eq:unify_1st}) to be
$N_\text{D} = N- \text{rank}(\Phi-w_\text{s}^\text{max} I_N )$,
where $w_\text{s}^\text{max}$ is the weight of the prevailing self-loop (see SM Sec. 2).
Using Eq.~(\ref{eq:unify_1st}) and the fact that $\Phi$ is a mixed matrix, we can
predict that $N_\text{D}$ remains unchanged if we exchange the
densities of any two types of self-loops (see SM Sec. 2), accounting for the symmetry of
$N_\text{D}$ for arbitrary types of self-loops. Due to the dominance of
$N_\text{D}$ by the self-loop with the highest density and the exchange-invariance
of $N_\text{D}$, the highest controllability with the lowest value of $N_\text{D}$ emerges
when distinct self-loops are of the same density.

\begin{figure}
\begin{center}
\epsfig{figure=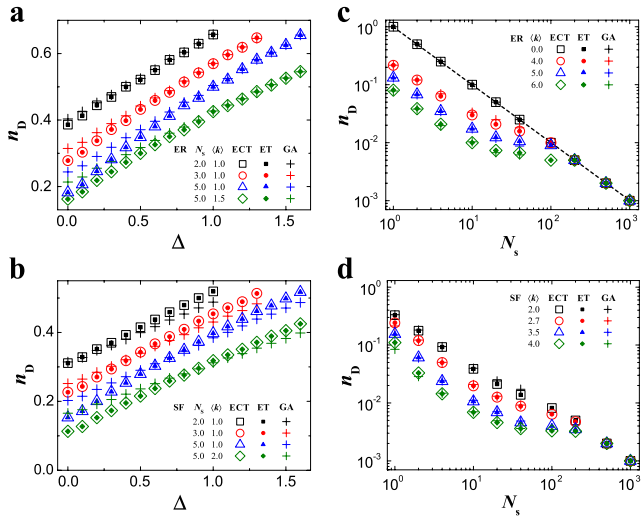,width=\linewidth}
\caption{{\bf a}-{\bf b}, $n_\text{D}$
as a function of the density heterogeneity of self-loops ($\Delta$) for ER ({\bf a}) and SF ({\bf b})
networks. {\bf c}-{\bf d}, $n_\text{D}$ as a function of the number of different types of self-loops
for ER ({\bf c}) and SF ({\bf d}) networks. The dotted line in ({\bf g}) is
$n_\text{D} = 1/ N_\text{s}$. The networks are described by structured matrix $A$ and their
sizes in ({\bf a})-({\bf d}) are 1000. The results from ECT and ET are averaged over 30 different realizations,
and those from GA are over 200 realizations. The notations are the same as Fig.~\ref{fig:1st_order}.
}
\label{fig:multi_dynamics}
\end{center}
\end{figure}

To validate the symmetry-induced highest controllability predicted by our theory, we quantify
the density heterogeneity of self-loops as follows:
\begin{equation}
\Delta \equiv \sum_{i=1}^{N_\text{s}}\left| \rho_\text{s}^{(i)}-\frac{1}{N_\text{s}} \right|,
\end{equation}
where $N_\text{s}$ is the number of different types of self-loops
(or the diversity of self-loops). Note that $\Delta=0$
if and only if all different types of self-loops have the same density, i.e., $\rho_\text{s}^{(1)}=
\rho_\text{s}^{(2)}=\cdots \rho_\text{s}^{(N_\text{s})}=\frac{1}{N_\text{s}}$,
and the larger value of $\Delta$ corresponds
to more diverse case. Figure~\ref{fig:multi_dynamics}a,b shows that $n_\text{D}$ monotonically
increases with $\Delta$ and the highest controllability (lowest $n_\text{D}$) arises at $\Delta=0$,
in exact agreement with our theoretical prediction.
Figure~\ref{fig:multi_dynamics}c,d display $n_\text{D}$ as a function of $N_\text{s}$. We see that
$n_\text{D}$ decreases as $N_\text{s}$ increases, suggesting that the diversity of individual dynamics
facilitates the control of a networked system. When $N_\text{s}=N$ (i.e., all the self-loops are independent),
$n_\text{D}=1/N$, which is also consistent with the prediction of SCT~\cite{LSB:2011,Cowan:2012}.\\

In some real networked systems, dynamic units are captured by high-order individual dynamics,
prompting us to check if the symmetry-induced highest controllability still holds
for higher-order individual dynamics.
The graph representation of dynamic units with 2nd-order dynamics
is illustrated in Fig.~\ref{fig:illus_motif}b. In this case,
the eigenvalues of the dynamic unit's state matrix $\left(
                                                      \begin{array}{cc}
                                                        0 & 1 \\
                                                        a_0 & a_1 \\
                                                      \end{array}
                                                    \right)$
play a dominant role in determining $N_\text{D}$.
For two different units as distinguished
by distinct ($a_0$ $a_1$) one can show that their state matrices almost
always have different eigenvalues, except for some pathological cases
of zero measure that occur when the parameters satisfy certain accidental constraints. The eigenvalues
of the state matrix of dynamic units take over the roles of self-loops in the 1st-order dynamics, accounting for the following
formulas for sparse networks
\begin{equation}
N_\text{D} = 2N- \min_i\bigg\{ \text{rank}(\Phi - \lambda^{(i)} I_{2N}) \bigg\},
\label{eq:2nd_nodal_formu}
\end{equation}
where $\lambda^{(i)}$ is either one of the two eigenvalues of
type-$i$ dynamic unit's state matrix. The formula implies that
$N_\text{D}$ is exclusively determined by the prevailing dynamic unit,
(see SM Sec. 2). The symmetry of $N_\text{D}$, i.e.,
exchanging the densities of any types of dynamic units, does not alter $N_\text{D}$
(see SM Sec. 2), and the emergence of highest controllability at the global symmetry point
can be similarly proved as we did in the case of 1st-order individual dynamics.

The 3rd-order individual dynamics are graphically characterized
by a dynamic unit composed of three nodes (Fig.~\ref{fig:illus_motif}c),
leading to a $3N\times 3N$ state matrix (Fig.~\ref{fig:illus_motif}c).
We can generalize Eq.~(\ref{eq:2nd_nodal_formu}) to
arbitrary order of individual dynamics:
\begin{equation}
N_\text{D} = dN- \min_i\bigg\{ \text{rank}(\Phi - \lambda^{(i)}_d I_{dN}) \bigg\},
\label{eq:unify_any_order}
\end{equation}
where $d$ is the order of the dynamic unit, $\lambda^{(i)}_d$
is any one of the $d$ eigenvalues of type-$i$ dynamic units and $I_{dN}$
is the identity matrix of dimension $dN$. In analogy with the simplified formula for the 1st-order
dynamics, insofar as a type of individual dynamics prevails in the system,
Eq~(\ref{eq:unify_any_order}) is reduced to
$N_\text{D} = dN- \text{rank}(\Phi - \lambda^\text{max}_d I_{dN})$,
where $\lambda^\text{max}_d$ is one of the eigenvalues of the prevailing dynamic unit's state matrix. Similar to the case of 1st-order individual dynamics,
the global symmetry of controllability
and the highest controllability occurs at the global symmetry point
can be proved for individual dynamics of any order and arbitrary
network topology (see SM Sec.2 and 3 for theoretical derivations and see SM Sec. 4 for numerical and analytical results of high-order individual dynamics).
In summary, we map individual dynamics
into dynamic units that can be integrated into the matrix representation of the system,
offering a general paradigm to explore the joint effect of individual dynamics and network
topology on the system's controllability. The paradigm leads to
a striking discovery: the universal symmetry of controllability as reflected by
the invariance of controllability with respect to exchanging the fractions of any two
different types of individual dynamics, and the emergence of highest controllability
at the global symmetry point. These findings generally hold for arbitrary
networks and individual dynamics of any order.
The symmetry-induced highest controllability has immediate implications for devising and
optimizing the control of complex systems by for example, perturbing individual dynamics
to approach the symmetry point without the need to adjust network structure.


The theoretical paradigm and tools developed here also allow us to address a number of
questions, answers to which could offer further insights into the control
of complex networked systems.
For example, we may consider the impact of general parameter
dependency (e.g., link weight similarity), instead of focusing on
self-loops or individual dynamics. Our preliminary results show that
introducing more identical link weights will not affect the network
controllability too much, unless the network is very dense and almost
all link weights are identical (see SM Sec.5). We still lack a comprehensive understanding
of the impact of parameter dependency on structural controllability
for arbitrary complex networks.
Moreover, at the present we are incapable of tackling general
nonlinear dynamical systems in the framework of ECT, which is extremely challenging for both
physicists and control theorists. Nevertheless, we hope our approach
could inspire further research interests towards achieving ultimate
control of complex networked systems.

\end{document}